\documentclass[12pt]{article}
\usepackage{amssymb}
\usepackage{amsfonts}
\usepackage{color}

\newcommand{\nc}{\newcommand}
\nc{\thref}[1]{Theorem~\ref{theo:#1}}
\nc{\selabel}[1]{\label{sect:#1}}
\nc{\seref}[1]{Section~\ref{sect:#1}}
\nc{\lelabel}[1]{\label{lemm:#1}}
\nc{\leref}[1]{Lemma~\ref{lemm:#1}}
\nc{\prlabel}[1]{\label{prop:#1}}
\nc{\prref}[1]{Proposition~\ref{prop:#1}}
\nc{\colabel}[1]{\label{coro:#1}}
\nc{\coref}[1]{Corollary~\ref{coro:#1}}
\nc{\exlabel}[1]{\label{exam:#1}}
\nc{\exref}[1]{Example~\ref{exam:#1}}
\nc{\delabel}[1]{\label{defi:#1}}
\nc{\deref}[1]{Definition~\ref{defi:#1}}
\nc{\eqlabel}[1]{\label{equa:#1}}
\nc{\relabel}[1]{\label{rema:#1}}
\nc{\reref}[1]{Lemma~\ref{rema:#1}}
\providecommand{\operatorname}[1]{\mathrm{#1}\,}
\nc{\Hom}{\operatorname{Hom}}
\nc{\Mor}{\operatorname{Mor}}
\nc{\Aut}{\operatorname{Aut}}
\nc{\Ann}{\operatorname{Ann}}
\nc{\Ker}{\operatorname{Ker}}
\nc{\Trace}{\operatorname{Trace}}
\nc{\Char}{\operatorname{Char}}
\nc{\Mod}{\operatorname{Mod}}
\nc{\End}{\operatorname{End}}
\nc{\Spec}{\operatorname{Spec}}
\nc{\Span}{\operatorname{Span}}
\nc{\sgn}{\operatorname{sgn}}
\nc{\Id}{\operatorname{Id}}
\nc{\Com}{\operatorname{Com}}
\nc{\rank}{\operatorname{rank}}

\let\:=\colon

\newtheorem{de}{Definition}[section]
\newtheorem{lm}[de]{Lemma}
\newtheorem{pr}[de]{Proposition}
\newtheorem{co}[de]{Corollary}
\newtheorem{re}[de]{Remark}
\newtheorem{res}[de]{Remarks}
\newtheorem{te}[de]{Theorem}
\newtheorem{ex}[de]{Example}
\newtheorem{exs}[de]{Examples}

\def\bex{\begin{ex}}
\def\eex{\end{ex}}
\def\bexs{\begin{exs}}
\def\eexs{\end{exs}}
\def\bl{\begin{lm}}
\def\el{\end{lm}}
\def\bc{\begin{co}}
\def\ec{\end{co}}
\def\bt{\begin{te}}
\def\et{\end{te}}
\def\bpr{\begin{pr}}
\def\epr{\end{pr}}
\def\br{\begin{re}}
\def\er{\end{re}}
\def\brs{\begin{res}}
\def\ers{\end{res}}
\def\bd{\begin{de}}
\def\ed{\end{de}}
\def\be{\begin{equation}}
\def\ee{\end{equation}}
\def\bea{\begin{eqnarray*}}
\def\eea{\end{eqnarray*}}
\def\bp{\begin{proof}}
\def\ep{\end{proof}}

\def\qed{\hfill\Box}
\let\:=\colon
\usepackage{amssymb}

\begin{document}
\begin{center}
\Large
% TITLE GOES HERE
Sharpness of the Finsler-Hadwiger inequality
\end{center}

\begin{flushright}
Cezar Lupu\\
Department of Mathematics-Informatics\\
University of Bucharest \\
Bucharest, Romania RO-010014\\
\verb+lupucezar@yahoo.com+

\vspace{2 mm}
Cosmin Pohoa\c t\u a \\
Tudor Vianu National College \\
Bucharest, Romania RO-010014\\
\verb+pohoata_cosmin2000@yahoo.com+
\end{flushright}

\it dedicated to the memory of the great professor, Alexandru
Lupa\c s \normalfont

\section{Introduction \& Preliminaries}\selabel{0}

The Hadwiger-Finsler inequality is known in literature of mathematics
 as a generalization of the %%@
following

\bt\label{t1}
In any triangle $ABC$ with the side lenghts $a, b, c$ and $S$ its area,
 the following inequality %%@
is valid
$$a^{2}+b^{2}+c^{2}\geq 4S\sqrt{3}.$$
\et This inequality is due to Weitzenbock, Math. Z, 137-146, 1919,
but this has also appeared at International Mathematical Olympiad
in 1961.
In [7.], one can find eleven proofs. In fact, in any triangle $ABC$ the
 following sequence of %%@
inequalities is valid:
$$a^2+b^2+c^2\geq ab+bc+ca\geq a\sqrt{bc}+b\sqrt{ca}+c\sqrt{ab}\geq
 3\sqrt [3] {a^2b^2c^2}\geq %%@
4S\sqrt{3}.$$

A stronger version is the one found by Finsler and Hadwiger in
1938, which states that ([2.]) \bt \label{t2_cezar}
In any triangle $ABC$ with the side lenghts $a, b, c$ and $S$ its area,
 the following inequality %%@
is valid
$$a^{2}+b^{2}+c^{2}\geq 4S\sqrt{3}+(a-b)^{2}+(b-c)^{2}+(c-a)^{2}.$$
\et In [8.] the first author of this note gave a simple proof only
by using AM-GM and the following inequality due to Mitrinovic: \bt
\label{t3}
In any triangle $ABC$ with the side lenghts $a, b, c$ and $s$ its
 semiperimeter and $R$ its %%@
circumradius, the following inequality holds
$$s\leq\frac{3\sqrt{3}}{2}R.$$
\et

This inequality also appears in [3.].

A nice inequality, sharper than Mitrinovic and equivalent to the first
 theorem is the following:
\bt \label{t4}
In any triangle $ABC$ with sides of lenghts $a, b, c$ and with inradius
 of $r$, circumradius of %%@
$R$ and $s$ its semiperimeter the following inequality holds
$$4R+r\geq s\sqrt{3}.$$

\et

In [4.], Wu gave a nice sharpness and a generalization of the
Finsler-Hadwiger inequality.

Now, we give an algebraic inequality due to I. Schur ([5.]), namely

\bt \label{t5}
For any positive real numbers $\displaystyle x, y, z$ and
 $t\in\mathbb{R}$
the following inequality holds

$$x^{t}(x-y)(x-z)+y^{t}(y-x)(y-z)+z^{t}(z-y)(z-x)\geq 0.$$

\et

The most common case is $\displaystyle t=1$, which has the following
 equivalent form:

$$x^3+y^3+z^3+3xyz\geq xy(x+y)+yz(y+z)+zx(z+x)$$
which is equivalent to
$$x^3+y^3+z^3+6xyz\geq (x+y+z)(xy+yz+zx).$$
Now, using the identity $\displaystyle
 x^3+y^3+z^3-3xyz=(x+y+z)(x^2+y^2+z^2-xy-yz-zx)$ one can %%@
easily deduce that
$$2(xy+yz+zx)-(x^2+y^2+z^2)\leq\frac{9xyz}{x+y+z}.(*)$$

Another interesting case is $\displaystyle t=2$. We have

$$x^4+y^4+z^4+xyz(x+y+z)\geq xy(x^2+y^2)+yz(y^2+z^2)+zx(z^2+x^2)$$
which is equivalent to
$$x^4+y^4+z^4+2xyz(x+y+z)\geq (x^2+y^2+z^2)(xy+yz+zx).(**)$$

Now, let's rewrite theorem 1.2. as $$2(ab+bc+ca)-(a^2+b^2+c^2)\geq
 4S\sqrt{3}.(***)$$
By squaring $(***)$ and using Heron formula we obtain
$$4\left(\sum_{cyc}ab\right)^{2}+\left(\sum_{cyc}a^{2}\right)^{2}-4\left(\sum_{cyc}ab\right)\left(%%@
\sum_{cyc}a^{2}\right)\geq 3(a+b+c)\prod (b+c-a)$$
which is equivalent to
$$6\sum_{cyc}{a^2b^2} + 4\sum_{cyc}{a^2bc} + \sum_{cyc}{a^4} -
 4\sum_{cyc}{ab(a+b)} \geq %%@
3(a+b+c)\prod (b+c-a).$$
By making some elementary calculations we get
$$6\sum_{cyc}{a^2b^2} + 4\sum_{cyc}{a^2bc} + \sum_{cyc}{a^4} -
 4\sum_{cyc}{ab(a+b)} \geq %%@
3(a+b+c)\left(\sum_{cyc}{ab(a+b)} - \sum_{cyc}{a^3} - 2abc\right).$$
We obtain the equivalent inequalities
$$\sum_{cyc}{a^4} +\sum_{cyc}{a^2bc} \geq \sum_{cyc}{ab(a^2+b^2)}$$

$$a^2(a-b)(a-c)+b^2(b-a)(b-c)+c^2(c-a)(c-b)\geq 0,$$
which is nothing else than Schur's inequality in the particular case
 $\displaystyle t=2$.
In what follows we will give another form of Schur's inequality. That
 is

\bt \label{t6} For any positive reals $\displaystyle m, n, p$, the
 following inequality holds
$$\frac{mn}{p}+\frac{np}{m}+\frac{pm}{n}+\frac{9mnp}{mn+np+pm}\geq
 2(m+n+p).$$
\et

\it Proof. \normalfont We denote $\displaystyle x=\frac{1}{m},
 y=\frac{1}{n}$ and $\displaystyle %%@
z=\frac{1}{p}$. We obtain the equivalent inequality
$$\displaystyle\frac{x}{yz}+\frac{y}{zx}+\frac{z}{xy}+\frac{9}{x+y+z}\geq\frac{2(xy+yz+zx)}{xyz}
\Leftrightarrow$$
$$\displaystyle
 2(xy+yz+zx)-(x^{2}+y^{2}+z^{2})\leq\frac{9xyz}{x+y+z},$$
which is $(*)$.

\section{Main results}\selabel{1}
In the previous section we stated a sequence of inequalities stronger
 than Weitzenbock inequality. %%@
In fact, one can prove that the following sequence of inequalities
 holds
$$a^2+b^2+c^2\geq ab+bc+ca\geq a\sqrt{bc}+b\sqrt{ca}+c\sqrt{ab}\geq
 3\sqrt [3] {a^2b^2c^2}\geq %%@
18Rr,$$
where $\displaystyle R$ is the circumradius and $\displaystyle r$ is
 the inradius of the triangle %%@
with sides of lenghts $\displaystyle a, b, c$. In this moment, one
 expects to have a stronger %%@
Finsler-Hadwiger inequality with $\displaystyle 18Rr$ instead of
 $\displaystyle 4S\sqrt{3}$. %%@
Unfortunately, the following inequality holds true
$$a^{2}+b^{2}+c^{2}\leq 18Rr+(a-b)^{2}+(b-c)^{2}+(c-a)^{2},$$
because it is equivalent to
$$2(ab+bc+ca)-(a^{2}+b^{2}+c^{2})\leq 18Rr=\frac{9abc}{a+b+c},$$
which is $(*)$ again. Now, we are ready to prove the first
refinement of the Finsler-Hadwiger inequality:

\bt\label{t1}
In any triangle $ABC$ with the side lenghts $a, b, c$ with $S$ its
 area, $\displaystyle R$ the %%@
circumradius and $\displaystyle r$ the inradius of the triangle
 $\displaystyle ABC$ the following %%@
inequality is valid

$$a^{2}+b^{2}+c^{2}\geq
 2S\sqrt{3}+2r(4R+r)+(a-b)^{2}+(b-c)^{2}+(c-a)^{2}.$$
\et

\it Proof. \normalfont We rewrite the inequality as
$$\displaystyle 2(ab+bc+ca)-(a^{2}+b^{2}+c^{2})\geq
 2S\sqrt{3}+2r(4R+r).$$
Since, $\displaystyle ab+bc+ca=s^{2}+r^{2}+4Rr$, it follows immediately
 that $\displaystyle %%@
a^{2}+b^{2}+c^{2}=2(s^{2}-r^{2}-4Rr)$.
The inequality is equivalent to
$$\displaystyle 16Rr+4r^{2}\geq 2S\sqrt{3}+2r(4R+r).$$
We finally obtain
$$\displaystyle  4R+r\geq s\sqrt{3},$$ which is exactly theorem 1.4.

$\qed$

The second refinement of the Finsler-Hadwiger inequality is the
following

\bt\label{t2}
In any triangle $ABC$ with the side lenghts $a, b, c$ with $S$ its
 area, $\displaystyle R$ the %%@
circumradius and $\displaystyle r$ the inradius of the triangle
 $\displaystyle ABC$ the following %%@
inequality is valid

$$\displaystyle a^{2}+b^{2}+c^{2}\geq %%@
4S\sqrt{3+\displaystyle\frac{4(R-2r)}{4R+r}}+(a-b)^{2}+(b-c)^{2}+(c-a)^{2}.$$

\et

\it Proof. \normalfont In theorem 1.6 we put $\displaystyle
 m=\frac{1}{2}(b+c-a), %%@
n=\frac{1}{2}(c+a-b)$ and $\displaystyle p=\frac{1}{2}(a+b-c)$. We get

$$\displaystyle\sum_{cyc}\frac{(b+c-a)(c+a-b)}{(a+b-c)}+
\frac{9(b+c-a)(c+a-b)(a+b-c)}{\displaystyle\sum_{cyc}(b+c-a)(c+a-b)}\geq
 2(a+b+c).$$
Since $\displaystyle ab+bc+ca=s^{2}+r^{2}+4Rr$ $(1)$ and
 $a^{2}+b^{2}+c^{2}=2(s^{2}-r^{2}-4Rr)$ %%@
$(2)$, we deduce
$$\displaystyle\sum_{cyc}(b+c-a)(c+a-b)=4r(4R+r).$$
On the other hand, by Heron's formula we have $\displaystyle
 (b+c-a)(c+a-b)(a+b-c)=8sr^{2}$, so %%@
our inequality is equivalent to
$$\displaystyle\sum_{cyc}\frac{(b+c-a)(c+a-b)}{(a+b-c)} +
 \frac{18sr}{4R+r} \geq 4s %%@
\Leftrightarrow$$
$$\displaystyle\sum_{cyc}\frac{(s-a)(s-b)}{(s-c)} + \frac{9sr}{4R+r}
 \geq 2s \Leftrightarrow \displaystyle\sum_{cyc}{(s-a)^2(s-b)^2} +
 \frac{9s^2r^3}{4R+r} \geq 2s^2r^2.$$
Now, according to the identity
$$\sum_{cyc}{(s-a)^2(s-b)^2} = \left(\sum_{cyc}{(s-a)(s-b)}\right)^2 -
 2s^2r^2,$$
we have
$$\displaystyle\left(\sum_{cyc}{(s-a)(s-b)}\right)^2-2s^2r^2 +
 \frac{9s^2r^3}{4R+r} \geq 2s^2r^2.$$
And since 
$$\displaystyle\sum_{cyc}{(s-a)(s-b)}=r(4R+r),$$
it follows that
$$\displaystyle r^2(4R+r)^2 + \frac{9s^2r^3}{4R+r}  \geq 4s^2r^2,$$
which rewrites as
$$ \left(\frac{4R+r}{s}\right)^{2} + \frac{9r}{4R+r} \geq 4. $$
From the identities mentioned in $(1)$ and $(2)$ we deduce that
$$ \frac{4R+r}{s}=\frac{2(ab+bc+ca)-(a^2+b^2+c^2)}{4S}. $$
The inequality rewrites as
$$ \left(\frac{2(ab+bc+ca)-(a^2+b^2+c^2)}{4S}\right)^2 \geq 4 -
 \frac{9r}{4R+r} \Leftrightarrow %%@
$$
$$ \left(\frac{(a^2+b^2+c^2) - \left( (a-b)^2+ (b-c)^2 + (c-a)^2
 \right)}{4S}\right)^2 \geq 3 + %%@
\frac{4(R-2r)}{4R+r} \Leftrightarrow $$
$$ \displaystyle a^{2}+b^{2}+c^{2}\geq %%@
4S\sqrt{3+\displaystyle\frac{4(R-2r)}{4R+r}}+(a-b)^{2}+(b-c)^{2}+(c-a)^{2}.$$
 $\qed$

{\bf{Remark.}} From Euler inequality, $\displaystyle R\geq 2r$, we
 obtain theorem 1.2.

\section{Applications}
In this section we illustrate some basic applications of the second
 refinement of Finsler-Hadwiger %%@
inequality. We begin with\\

{\bf{Application 1.}} \it In any triangle $\displaystyle ABC$ with the
 sides of lenghts %%@
$\displaystyle a, b, c$ the following inequality holds\normalfont
$$\displaystyle\frac{1}{b+c-a}+\frac{1}{c+a-b}+\frac{1}{c+a-b}\geq\frac{1}{2r}\sqrt{4-\frac{9r}{4%%@
%%@
R+r}}.$$

\it Solution.\normalfont From $$\displaystyle
 (b+c-a)(c+a-b)(a+b-c)=4r(4R+r),$$ it is quite easy to observe that
$$\displaystyle\frac{1}{b+c-a}+\frac{1}{c+a-b}+\frac{1}{a+b-c}=\frac{4R+r}{2sr}.$$
Now, applying the inequality $$\displaystyle
 \left(\frac{4R+r}{s}\right)^{2}+\frac{9r}{4R+r}\geq %%@
4,$$ we get
$$\displaystyle\left(\frac{1}{b+c-a}+\frac{1}{c+a-b}+\frac{1}{a+b-c}\right)^{2}=\frac{1}{4r^{2}}
 %%@
\left(\frac{4R+r}{s}\right)^{2}\geq\frac{1}{4r^{2}}\left(4-\frac{9r}{4R+r}\right).$$
The given inequality follows immediately. $\qed$

{\bf{Application 2.}} \it In any triangle $\displaystyle ABC$ with the
 sides of lenghts %%@
$\displaystyle a, b, c$ the following inequality holds\normalfont
$$\displaystyle\frac{1}{a(b+c-a)}+\frac{1}{b(c+a-b)}+\frac{1}{c(a+b-c)}\geq\frac{r}{8R}\left(5-
 %%@
\frac{9r}{4R+r}\right).$$

\it Solution.\normalfont  From the following identity
$$\displaystyle\sum_{cyc}\frac{(s-a)(s-b)}{c}=\frac{r(s^{2}+(4R+r)^{2})}{4sR}=\frac{S}{4R}\left(1+%%@
\left(\frac{4R+r}{p}\right)^{2}\right).$$
Using the inequality
$$\displaystyle\left(\frac{4R+r}{s}\right)^{2}+\frac{9r}{4R+r}\geq 4,$$
we have
$$\displaystyle\sum_{cyc}\frac{(s-a)(s-b)}{c}\geq\frac{S}{4R}\left(5-\frac{9r}{4R+r}\right).$$
In this moment, the problem follows easily. $\qed$

{\bf{Application 3.}} \it In any triangle $\displaystyle ABC$ with the
 sides of lenghts %%@
$\displaystyle a, b, c$ the following inequality holds\normalfont
$$\displaystyle\frac{1}{(b+c-a)^{2}}+\frac{1}{(c+a-b)^{2}}+\frac{1}{(a+b-c)^{2}}\geq
 %%@
\frac{1}{r^{2}}\left(\frac{1}{2}-\frac{9r}{4(4R+r)}\right).$$\\

\it Solution.\normalfont
From $\displaystyle (b+c-a)(c+a-b)(a+b-c)=4r(4R+r),$ it follows that
$$(b+c-a)^{2}+(c+a-b)^{2}+(a+b-c)^{2}=4(s^2-2r^2-8Rr)$$
and
$$(b+c-a)^{2}(c+a-b)^{2}+(a+b-c)^{2}(c+a-b)^{2}+(b+c-a)^{2}(a+b-c)^{2}=4r^{2}\left((4R+r)^{2}-2s^{2}\right).$$
We get
$$\displaystyle\frac{1}{(b+c-a)^{2}}+\frac{1}{(c+a-b)^{2}}+\frac{1}{(a+b-c)^{2}}=\frac{1}{4}\left(\frac{(4R+r)^{2}}{s^{2}r^{2}}-\frac{2}{r^{2}}\right).$$
Now, applying the inequality
$$\displaystyle\left(\frac{4R+r}{s}\right)^{2}+\frac{9r}{4R+r}\geq 4,$$
we have
$$\displaystyle\frac{1}{(b+c-a)^{2}}+\frac{1}{(c+a-b)^{2}}+\frac{1}{(a+b-c)^{2}}\geq\frac{1}{4r^{2}}\left(2-\frac{9r}{4R+r}\right)=\frac{1}{r^{2}}\left(\frac{1}{2}-\frac{9r}{4(4R+r)}\right).$$
 $\qed$

{\bf{Application 4.}} \it In any triangle $\displaystyle ABC$ with the
 sides of lenghts %%@
$\displaystyle a, b, c$ the following inequality holds\normalfont
$$\displaystyle\frac{a^2}{b+c-a}+\frac{b^2}{c+a-b}+\frac{c^2}{a+b-c}\geq
 %%@
3R\sqrt{4-\frac{9r}{4R+r}}.$$\\

\it Solution.\normalfont
Without loss of generality, we assume that $\displaystyle a\leq b\leq
 c$. It follows quite easily %%@
that $\displaystyle a^2\leq b^2\leq c^2$ and %%@
$\displaystyle\frac{1}{b+c-a}\leq\frac{1}{c+a-b}\leq\frac{1}{a+b-c}$.
 Applying Chebyshev's %%@
inequality, we have
$$\displaystyle\frac{a^2}{b+c-a}+\frac{b^2}{c+a-b}+\frac{c^2}{a+b-c}\geq\left(\frac{a^2+b^2+c^2}{3%%@
}
\right)\left(\frac{1}{b+c-a}+\frac{1}{c+a-b}+\frac{1}{c+a-b}\right).$$
Now, the first application and the inequality $\displaystyle
 a^2+b^2+c^2\geq 18Rr$ solves the %%@
problem. $\qed$

{\bf{Application 5.}} \it In any triangle $\displaystyle ABC$ with the
 sides of lenghts %%@
$\displaystyle a, b, c$ and with the exradii $\displaystyle r_{a},
 r_{b}, r_{c}$ corresponding to the triangle $ABC$, the following inequality
 holds\normalfont
$$\displaystyle\frac{a}{r_{a}}+\frac{b}{r_{b}}+\frac{c}{r_{c}}\geq %%@
2\sqrt{3+\frac{4(R-2r)}{4R+r}}.$$

\it Solution.\normalfont
From the well-known relations $\displaystyle r_{a}=\frac{S}{s-a}$ and
 %%@
the analogues, the inequality is equivalent to
$$\displaystyle\frac{a}{r_{a}}+\frac{b}{r_{b}}+\frac{c}{r_{c}}=\frac{2(ab+bc+ca)-(a^2+b^2+c^2)}{2S%%@
}\geq 2S\sqrt{3+\frac{4(R-2r)}{4R+r}}.$$
The last inequality follows from theorem 2.2 immediately. $\qed$

{\bf{Application 6.}} \it In any triangle $\displaystyle ABC$ with the
 sides of lenghts %%@
$\displaystyle a, b, c$ and with the exradii $\displaystyle r_{a},
 r_{b}, r_{c}$ corresponding to the triangle $ABC$ and with $\displaystyle
 h_{a}, h_{b}, h_{c}$ be the altitudes of the  triangle $\displaystyle
 ABC$, the following inequality holds\normalfont

$$\displaystyle\frac{1}{h_{a}r_{a}}+\frac{1}{h_{b}r_{b}}+\frac{1}{h_{c}r_{c}}\geq\frac{1}{S}\sqrt{3+\frac{4(R-2r)}{4R+r}}.$$

\it Solution.\normalfont From the well-known relations in triangle
$ABC$, $\displaystyle h_{a}=\frac{2S}{a},\displaystyle
r_{a}=\frac{S}{s-a}$ we have
$\displaystyle\frac{1}{h_{a}r_{a}}=\frac{a(s-a)}{2S^{2}}$. Doing
the same thing for the analogues and adding them up we get
$$\displaystyle\frac{1}{h_{a}r_{a}}+\frac{1}{h_{b}r_{b}}+\frac{1}{h_{c}r_{c}}=\frac{1}{2S^{2}}\left(a(s-a)+b(s-b)+c(s-c)\right).$$
On the other hand by using theorem 2.2 in the form
$$\displaystyle a(s-a)+b(s-b)+c(s-c)\geq
 2S\sqrt{3+\frac{4(R-2r)}{4R+r}}$$
we obtain the desired inequality. $\qed$

{\bf{Application 7.}} \it In any triangle $\displaystyle ABC$ with the
 sides of lenghts $\displaystyle a, b, c$ the following inequality holds
 true

$$\displaystyle\tan\frac{A}{2}+\tan\frac{B}{2}+\tan\frac{C}{2}\geq\sqrt{3+\frac{4(R-2r)}{4R+r}}.$$

\it Solution.\normalfont
From the cosine law we get $\displaystyle a^{2}=b^{2}+c^{2}-2bc\cos A$.
Since $S=\frac{1}{2}bc\sin A$ it follows that
$$\displaystyle a^{2}=(b-c)^{2}+4S\cdot\frac{1-\cos A}{\sin A}.$$
On the other hand by the trigonometric formulae $\displaystyle 1-\cos
 A=2\sin^{2}\frac{A}{2}$ and $\displaystyle\sin
 A=2\sin\frac{A}{2}\cos\frac{A}{2}$ we get
$$\displaystyle a^{2}=(b-c)^{2}+4S\tan\frac{A}{2}.$$
Doing the same for all sides of the triangle $ABC$ and adding up we
 obtain
$$a^{2}+b^{2}+c^{2}=(a-b)^{2}+(b-c)^{2}+(c-a)^{2}+4S\left(\tan\frac{A}{2}+\tan\frac{B}{2}+\tan\frac{C}{2}\right).$$
Now, by theorem 2.2 the inequality follows. $\qed$

{\bf{Application 8.}} \it In any triangle $\displaystyle ABC$ with the
 sides of lenghts %%@
$\displaystyle a, b, c$ and with the exradii $\displaystyle r_{a},
 r_{b}, r_{c}$ corresponding to the triangle $ABC$, the following inequality
 holds\normalfont
$$\displaystyle\frac{r_{a}}{a}+\frac{r_{b}}{b}+\frac{r_{c}}{c}\geq
 \frac{s(5R-r)}{R(4R+r)}.$$

\it Solution.\normalfont It is well-known that the following identity
 is valid in any triangle ABC
$$\displaystyle\frac{r_{a}}{a}+\frac{r_{b}}{b}+\frac{r_{c}}{c}=\frac{(4R+r)^2+s^2}{4Rs}.$$
So, the inequality rewrites as
$$\displaystyle\frac{(4R+r)^2}{s^2} + 1 \geq \frac{4(5R-r)}{4R+r},$$
which is equivalent with
$$\displaystyle\left(\frac{4R+r}{s}\right)^{2}+\frac{9r}{4R+r}\geq 4.$$
 $\qed$

\paragraph*{Acknowledgment.}
The authors would like to thank to Nicolae Constantinescu, from
 University of Craiova and to %%@
Marius Ghergu, from the Institute of Mathematics of the Romanian
 Academy for useful suggestions. This paper has been completed while the first
 author participated in the summer school on \it Critical Point theory
 and its applications \normalfont organized in Cluj-Napoca city. We are
 kindly grateful to professors Vicen\c tiu R\u adulescu from the
 Institute of Mathematics of the Romanian Academy and to Csaba Varga from
 Babe\c s-Bolyai University, Cluj-Napoca.  \\

\end{document}